\documentclass[10pt,a4paper,twocolumn]{article}

\usepackage[T1]{fontenc}
\usepackage[utf8]{inputenc}
\usepackage[english]{babel}
\usepackage[margin=2cm]{geometry}
\usepackage{amsmath,amssymb,amsthm}
\usepackage{graphicx}
\usepackage{xcolor}
\usepackage[hidelinks]{hyperref}
\usepackage{generic}
\usepackage[table]{xcolor}
\usepackage{cite}

\usepackage{tcolorbox}

\usepackage{amsthm}

\usepackage{algorithmic}
\usepackage{graphicx}
\usepackage{algorithm,algorithmic}
\usepackage{hyperref}
\hypersetup{hidelinks}
\usepackage{textcomp}

\def\BibTeX{{\rm B\kern-.05em{\sc i\kern-.025em b}\kern-.08em
    T\kern-.1667em\lower.7ex\hbox{E}\kern-.125emX}}

\date{}

\begin{document}
\title{Twofold Control: Simultaneous Optimization of Initial Conditions and Control Inputs}
\author{Martin Lazar,  Mara Vlašić
\thanks{This work was supported by the Croatian Science Foundation under the project number HRZZ-IP-2022-10-5191. This work was funded by the European Union - NextGenerationEU through the institutional project TRIUMPh at the University of Dubrovnik.}
\thanks{The authors are with the Center for Control Theory and Applied Mathematics, Faculty of Electrical Engineering and Applied Computing, University of Dubrovnik, Dubrovnik, Croatia (emails: martin.lazar@unidu.hr, mara.vlasic@unidu.hr).}}
\maketitle

\begin{abstract}
The paper introduces a novel control notion, twofold controllability, which consists of the joint optimisation of the initial state and the distributed control input in a controlled dynamical system, with the aim of reaching a prescribed target state within a given time interval. We show that the resulting problem is well-posed, derive the corresponding optimality conditions, and provide an explicit formulation of its solution. The proposed approach ensures the reachability of any target state independently of the Kalman rank condition. We show that incorporating the additional control component enlarges the set of states reachable with a fixed control cost. At the same time, the average energy required to reach target vectors on the unit sphere decreases. Furthermore, we analyse the asymptotic behaviour as the cost of controlling the system through the initial state tends to infinity. Particular attention is paid to the case in which the controllability Gramian is singular and the classical notion of controllability does not hold. The theoretical results are illustrated by numerical examples, including ones arising from the semi-discretisation of the heat and wave equations with boundary control.
\end{abstract}

\begin{KeyWords}
    Twofold control, Controllability Gramian, Controllability measures, Optimal cost
\end{KeyWords}

\section{Introduction}\label{sec:introduction}

Controllability is a fundamental concept in control theory and concerns the ability to steer a dynamical system to a desired state by means of an appropriate control input. In the standard formulation, the initial state is prescribed and the control action is applied exclusively through a time-dependent control input. The associated control cost and the set of reachable states are closely related to the controllability Gramian, which plays a central role in optimal control and in the quantitative analysis of controllability.

In many practical applications, however, the initial configuration of a system can also be adjusted before the control process begins. This motivates a more general control framework in which the initial state and the control input are optimized simultaneously. We refer to this setting as \emph{twofold control}. The main idea is to distribute the available control effort between preparing a suitable initial state and applying the subsequent control input. This additional degree of freedom may substantially reduce the required control energy and enlarge the set of reachable states.

To describe this approach in a finite-dimensional linear setting, we introduce the twofold controllability Gramian, which incorporates the contributions of both the optimally selected initial state and the control input. This extended Gramian naturally characterizes the reachable set under a common energy constraint and provides a direct generalization of the classical controllability Gramian. It also measures the average control cost required to reach a target state on the unit sphere. Building on these properties, we investigate the geometric, energetic, and optimization aspects of the twofold controllability problem and compare our findings with those obtained in the classical control setting.

The main contributions of this paper are as follows:
\begin{itemize}
    \item We formulate the twofold optimal control problem for  a finite-dimensional, linear system, in which the initial state and the control input are optimised simultaneously, derive the corresponding optimality system, and introduce an extended controllability Gramian that characterises the minimum control cost and the set of reachable terminal states.
     \item We analyse the geometric and energetic properties of the twofold control setting,  and compare them with the  classical one.
    \item We analyse the asymptotic behaviour of the optimal solution as the cost of initial-state preparation tends to infinity, considering both positive definite and singular controllability Gramians.
          \item The theoretical results are illustrated by numerical experiments for several finite-dimensional systems, including semi-discrete models of the one-dimensional heat and wave equations with boundary control.
\end{itemize}

The remainder of the paper is organised as follows. Section II-A introduces the twofold control problem, while Section II-B derives the corresponding optimality conditions. Section II-C investigates geometric and energetic properties of the twofold controllability, highlighting the advantages of the proposed framework over classical controllability measures. Section III  investigates the limiting behaviour as the penalty on the initial state tends to infinity, distinguishing between controllable and uncontrollable systems in the classical setting. The theoretical findings are verified through numerical experiments in Section IV. Finally, Section V concludes the paper.

\section{Twofold Control}

\subsection{Set-up of the Problem}\label{sec:setup}

We analyse finite-dimensional linear systems and  focus on optimal control problems in which both the initial state and the distributed control input are optimised simultaneously. The motivation is to reduce the control cost and ensure controllability independently of the Kalman rank condition. We begin with the problem setup and the formulation of the joint optimisation.

We consider the linear control system
\begin{align}\label{eq:system}
    \begin{aligned}
        y'(t) &= A y(t) + B u(t), \quad t \in [0,T], \\
        y(0) &= y^0,
    \end{aligned}    
\end{align}
where
\begin{itemize}
    \item $y(t) \in \mathbb{R}^n$ is the state vector,
    \item $u(t) \in \mathbb{R}^m$ is the distributed control input,
    \item $A \in M_n(\mathbb{R})$ is the system dynamics matrix,
    \item and $B \in M_{n,m}(\mathbb{R})$ is the control operator.
\end{itemize}
In the standard optimal control problem, one controls the system via the  control input $u$. However, before the control horizon begins (i.e., before the start of an experiment), one may prepare the system in advance by choosing some suitable initial state.

This leads to a new approach in which we optimise both the initial state and a distributed control input, aiming for greater flexibility and potentially lower cost. We refer to this notation as the {\it twofold control}.

Its goal is to steer the system to a given target final state \(y^T\) in a prescribed time horizon $T$, while minimising the cost functional
\begin{align}\label{eq:cost_ab}
    J(y^0,u) =  \alpha \| y^0 \|_{ \mathbb{R}^n }^2 + \beta \| u \|_{ L^2( 0,T; \mathbb{R}^m) }^2 , \quad \alpha, \beta > 0.
\end{align}
Unlike standard control problems, we optimise with respect to both the initial state \(y^0\) and the control input \(u(t)\).

The parameters \(\alpha\) and \(\beta\) are weights in the cost functional: \(\alpha\) penalises the size of the chosen initial state \(y^0\) and measures the cost of adjusting the initial state before the start of the control process. The second parameter \(\beta\) penalises the energy of the distributed control input \(u\). By tuning these parameters, we decide whether to prioritise keeping the initial configuration small or reducing the control effort.

\if 0 %beta=1 version
Of course,  the problem of minimising the cost functional $J$ does not change with rescaling the weights  \(\alpha\) and \(\beta\) by the same constant. For this reason, without losing generality, in the sequel, we fix the second parameter to $\beta=1$ and denote the corresponding functional by $J_\alpha$ (in order to stress its dependence on $\alpha$, what will be of important in the subsequent sections). 

Consequently, we formulate our problem as follows. 
\begin{align}\label{eq:cost_al}
    J_\alpha(y^0,u) = \frac{1}{2} \left( \alpha \| y^0 \|_{ \mathbb{R}^n }^2 + \beta \| u \|_{ L^2( 0,T; \mathbb{R}^m) }^2 \right), \quad \alpha, \beta > 0.
\end{align}
\fi 

To precisely formulate the optimisation problem in a functional-analytic setting, we introduce the control space \(X:=\mathbb{R}^n\times L^2(0,T;\mathbb{R}^m)\), defined as a Hilbert space equipped with the weighted inner product
\begin{align}\label{eq:innerproduct}
    \langle (y^0,u),(z^0,v)\rangle_X = \alpha\langle y^0,z^0\rangle_{\mathbb R^n} + \beta\langle u,v\rangle_{L^2(0,T;\mathbb R^m)},
\end{align}
and the corresponding norm
\begin{align}\label{eq:norm}
    \| (y^0, u) \|_X^2 = \alpha \| y^0 \|_{ \mathbb{R}^n }^2 + \beta \| u \|_{ L^2(0, T; \mathbb{R}^m) }^2.
\end{align}
The terminal state corresponding to a control pair \((y^0,u)\in X\) is given by
\begin{align}\label{eq:L(y,u)}
    L(y^0, u) = e^{T A} y^0 + \int_0^T e^{(T - t)A} B u(t) \, dt,
\end{align}
where \(L:X\to\mathbb R^n\) denotes the control-to-state operator. 

Consequently, we formulate our problem as follows. 

\begin{tcolorbox}[
    colback=white,
    colframe=black,
    boxrule=0.2pt,
    arc=0pt
]
{\sl Optimal twofold control problem:}
\begin{align}
\label{eq:min_J}
 \min_{(y^0, u) \in X}  J(y^0,u) 
\end{align}
subject to
\begin{align}\label{eq:Lu=xT}
    L(y^0, u) = y^T.
\end{align}
\end{tcolorbox}

\begin{nap}\label{nap-J_X}
    As \(J(y^0,u)= \|(y^0, u)\|_X^2\), the minimisation problem \eqref{eq:min_J} can be equivalently written as \(\min_{(y^0, u) \in X}  \|(y^0, u)\|_X^2\). 
\end{nap}

The well-posedness of the introduced problem is ensured by the following result. 

\begin{teo}\label{well_pos}
    The optimal twofold control problem admits a unique minimizer for every choice of the operators $ A \in M_n(\mathbb{R}), B \in M_{n,m}(\mathbb{R})$, every target  $y^T \in \mathbb{R}^n$, and every  weights $\alpha,\beta \in \mathbb{R}^+$. 
\end{teo}

\begin{proof}
    We introduce the set of admissible controls defined by
    \begin{equation*}
        \label{eq:A_adm}
        \mathcal A := \left\{ (y^0,u)\in X \;\middle|\; L(y^0,u)=y^T \right\},
    \end{equation*}
        that is, the set of all initial states and distributed controls steering the system \eqref{eq:system} to the prescribed target state \(y^T\).
    
    The admissibility constraint \eqref{eq:Lu=xT} can be incorporated into the objective functional by introducing the indicator function \(I_{\mathcal A}:X\longrightarrow\mathbb R\cup\{+\infty\}\), defined by
    \[
        I_{\mathcal A}(y^0,u)= 
        \begin{cases}
            0, & (y^0,u)\in\mathcal A,\\[1mm] +\infty, & (y^0,u)\notin\mathcal A.
        \end{cases}
    \]
    
    Consequently, the original constrained optimisation problem
    \[
        \min_{(y^0,u)\in\mathcal A}J(y^0,u)
    \]
    can be written equivalently in the unconstrained form
    \[
        \min_{(y^0,u)\in X} F(y^0,u),
    \]
    where
    \begin{align}\label{eq:F=J+IA}
        F(y^0,u) := J(y^0,u) + I_{\mathcal A}(y^0,u).
    \end{align}
    Since the admissible set \(\mathcal A\) is a non-empty, closed, and convex subset of the Hilbert space \(X\), its indicator function \(I_{\mathcal A}\) is convex and lower semicontinuous \cite[\textit{Chapter 2}]{JP15}. Consequently, the functional \(F(y^0,u)\) is proper, coercive, lower semicontinuous, and strictly convex. Therefore, by \cite[Theorem 2.19]{JP15}, the optimisation problem admits a unique minimiser.
\end{proof}

In the classical control setting, the linear control system \eqref{eq:system} is controllable, if the pair $(A,B)$ satisfies the Kalman rank condition \cite[Theorem 3]{ES13},
\[
    \operatorname{rank} \left[ B, \, AB, \, A^2B, \ldots, A^{n-1} B \right] = n.
\]

In contrast, the twofold control problem considered in this paper does not require the Kalman rank condition. As shown above,  the proposed optimisation problem is well-posed for any pair \((A,B)\), independently of whether the corresponding classical control system is controllable. In particular, the simultaneous optimisation of the initial state and the control input allows any prescribed final state to be reached within a given time horizon, even when the Kalman rank condition is not satisfied.

\subsection{Optimality Conditions for Twofold Control}\label{sec-optim}

Having established the well-posedness of the problem, we seek for the expressions determining the optimal control pair. To this end, we recall the previously defined control-to-state operator \(L\) and introduce its adjoint \(L^\star\). The operator \(L\) maps a pair consisting of the initial state and the distributed control input to the corresponding final state, whereas \(L^\star\) maps a final-state of the dual variable to the  control space $X$. This operator formulation will be used to characterise the optimal solution and derive the corresponding optimality conditions.

Motivated by the Hilbert uniqueness method (HUM), we seek a feasible pair of the form 
\begin{align}\label{eq:(y0,u)=L*pT}
    (y^0, u) = L^\star p^T
\end{align}
where \(L^\star\) stands for the adjoint operator  defined by
\begin{align}\label{eq:def_Lstar}
    \langle L(y^0, u), p^T \rangle_{\mathbb{R}^n} = \langle (y^0, u), L^\star p^T \rangle_X,
\end{align}
for every \((y^0, u) \in X\) and every \( p^T \in \mathbb{R}^n\) \cite[Prop. 4.28]{JvN23}.

Using the explicit representation of the controllability operator \(L\) in \eqref{eq:L(y,u)}, and transferring the operators to the second argument of the inner product on the left side of \eqref{eq:def_Lstar}, we obtain
\begin{align*}
    \langle L(y^0,u) ,p^T \rangle_{\mathbb R^n} 
    &= \langle y^0,e^{TA^\star}p^T\rangle_{\mathbb R^n} \\
    &\quad + \int_0^T \langle u(t),B^\star e^{(T-t)A^\star}p^T\rangle_{\mathbb R^m}\,dt.
\end{align*}
Comparison of this expression with the weighted inner product on \(X\) for \(\langle (y^0,u),L^\star p^T\rangle_X\) yields the explicit formula for the adjoint operator
\begin{align}\label{eq:Lstar}
    L^\star p^T = \left( \frac{1}{\alpha} e^{T A^\star} p^T, \frac{1}{\beta} B^\star e^{(T-t)A^\star} p^T \right).
\end{align}
Consequently, the control pair given by \eqref{eq:(y0,u)=L*pT} satisfies
\begin{align}\label{eq:x0u}
    y^0= \frac{1}{\alpha} e^{T A^\star} p^T, \, u(t)= \frac{1}{\beta} B^\star e^{(T-t)A^\star} p^T, \, t \in [0, T]
\end{align}
and the corresponding system  final state $L(y^0, u)$ can be expressed in a compact form using the Gramian operator. 

By combining \eqref{eq:x0u} and \eqref{eq:Lu=xT}, we obtain the 
\begin{align}\label{eq:xT=MpT}
     Mp^T=y^T ,
\end{align}
where $M$ denotes the {\it twofold} (or { \it extended}) { \it controllability Gramian}  defined by
\begin{align}\label{eq:mca_M}
    M: = L L^\star =\frac{1}{\alpha}E + \frac{1}{\beta} G,
\end{align}
with
\[
   E = e^{TA}e^{TA^\star}, \quad  G = \int_0^T e^{(T-t)A}BB^\star e^{(T-t)A^\star} \, dt.
\]

The latter term corresponds to the distributed control input and coincides with the classical controllability Gramian \(G\), while the first one, \(E\), quantifies the contribution of the optimisable initial state.

As \(E\) and \(G\) are positive definite and semidefinite, respectively, for every \(\alpha,\beta>0\), the operator \(M\) is positive definite and therefore invertible. Consequently, for every target state \(y^T\in\mathbb R^n\), there exists a unique adjoint datum \(\tilde p^T\) satisfying \eqref{eq:xT=MpT}. 

The next result shows that this adjoint datum also provides  the optimal control pair, uniquely determined by \eqref{eq:x0u}.

\begin{teo}\label{thm-optim}
    Let \(y^T\in\mathbb R^n\) be an arbitrary target final state and let \(\tilde p^T=M^{-1}y^T\). Then \((\tilde y^0,\tilde u)=L^\star \tilde p^T\) is the unique solution of the optimal twofold control problem
    \[
        \min_{\substack{
        (y^0,u) \in \mathcal A
    }}J(y^0,u).
    \]
    The associated, minimal  control cost is 
    \begin{align}\label{opt_cost}
        \|(\tilde y^0,\tilde u)\|_X^2 = (y^T)^\top M^{-1} y^T
    \end{align}
\end{teo}

\begin{proof}
    First note that by the invertibility of the twofold controllability Gramian, the state $\tilde p^T$ is unique and well defined. Furthermore, the corresponding control pair  is feasible, i.e.,
    \[
        L(\tilde y^0,\tilde u)=LL^\star \tilde p^T=M \tilde p^T=y^T.
    \]
    
    Let \((z^0,v)\) be any other feasible pair and define \(h=(z^0-\tilde y^0,v-\tilde u)\). Then \(Lh=0\). Using the definition of the adjoint operator \eqref{eq:def_Lstar}, we obtain
    \[
        \langle (\tilde y^0, \tilde u),h\rangle_X = \langle L^\star \tilde p^T,h\rangle_X = \langle \tilde p^T,Lh\rangle_{\mathbb R^n} = 0.
    \]
    Therefore,
    \[
        \|(z^0,v)\|_X^2 = \|(\tilde y^0,\tilde u)\|_X^2+\|h\|_X^2 \ge \|(\tilde y^0,\tilde u)\|_X^2,
    \]
    which, according to Remark \ref{nap-J_X}, implies \(J(z^0,v)\ge J(\tilde y^0,\tilde u)\).
    
    Moreover, equality holds if and only if \(h=0\), proving that \((\tilde y^0,\tilde u)\) is the unique minimiser.

    In order to show \eqref{opt_cost}, note that the energy of the optimal pair is
    \begin{align*}
       \|(\tilde y^0,\tilde u)\|_X^2 &= \|L^\star \tilde p^T\|_X^2 \\
        &= \langle LL^\star \tilde p^T,\tilde p^T\rangle_{\mathbb{R}^n} = (\tilde p^T)^\top M \tilde p^T.
    \end{align*}
    
    Substituting \(\tilde p^T = M^{-1}y^T\) yields the result. 
\end{proof}

Positive definiteness of the twofold controllability Gramian guarantees the existence of a unique optimal adjoint state \(\tilde p^T\), whereas Theorem \ref{thm-optim} provides explicit formulae for the corresponding optimal pair \((\tilde y^0,\tilde u)\).

\begin{nap}
    The optimality system derived above is not specific to the HUM approach adopted in this paper. The same optimality conditions can also be obtained by the method of Lagrange multipliers or, equivalently, by the Pontryagin Maximum Principle after reformulating the problem as a standard optimal control problem; see, e.g., \cite{KK08,ET23}. This confirms the consistency of the proposed formulation and shows that the derived optimality conditions are independent of the particular derivation technique.
\end{nap}

\subsection{Geometric and Energetic Properties of  Twofold Controllability  }

The controllability measures based on the determinant and the trace of the inverse of the controllability Gramian have been studied in \cite{SM21,SC21,SS16,FL16,LA21} for the classical control problem, in which the initial state is fixed and controllability is achieved solely through the distributed control input. These measures have natural geometric and energetic interpretations: the determinant is related to the volume of the reachable set under a unit-energy constraint, whereas the trace of the inverse quantifies the average minimum energy required to reach terminal states uniformly distributed on the unit sphere.

In the twofold control framework considered in this article, the initial state is not prescribed but is instead optimised jointly with the control input. In order to investigate the reachability properties of this approach, we subject both quantities to the common energy constraint
\begin{align}\label{eq:energy_constraint}
   \| (y^0, u) \|_X^2 = \alpha\|y^0\|_{\mathbb{R}^n}^2 + \beta\int_0^T \|u(t)\|_{\mathbb{R}^m}^2\,dt \leq 1,
\end{align}
where $\alpha,\beta>0$ determine the relative cost of initial-state preparation and distributed control effort.

The corresponding reachable set associated with the cost constraint \eqref{eq:energy_constraint} is defined by
\begin{align*}
    \mathcal{E}_{\mathrm{ext}} = \left\{ y^T \in \mathbb{R}^n : 
    %\right. & \exists (y^0,u) \in X, \\    &\left.
        y^T = L(y^0,u), \; \|(y^0,u)\|_X^2 \leq 1 \right\}
\end{align*}

According to Theorem \ref{thm-optim}, the minimal cost required for reaching any prescribed target \(y^T\) is obtained by the HUM control of the form \((\tilde y^0,\tilde u)=L^\star M^{-1}y^T\), and it is enough  to restrict the definition of the reachable set \(\mathcal{E}_{\mathrm{ext}}\) to HUM controls. According to the same result, for such controls, the cost is equal to \((y^T)^\top M^{-1} y^T\). 

Consequently, the reachable set has the ellipsoidal representation
\begin{align}\label{eq:Eext}
    \mathcal{E}_{\mathrm{ext}} = \left\{ y^T \in \mathbb{R}^n : (y^T)^\top M^{-1} y^T \leq 1 \right\}.
\end{align}

Since the volume of the ellipsoid \eqref{eq:Eext} is proportional to \(\sqrt{\det\bigl(M\bigr)}\), as shown in \cite[Corollary 12.15]{NV21}, the determinant of the extended Gramian measures the reachability potential in the twofold control framework. 

Apart from analysing the volume of the  reachable set, one is also interested in lowering the control cost. As shown above, the minimal cost of reaching any target $y^T$  equals \((y^T)^\top M^{-1} y^T\). On the other hand,  for any symmetric matrix \(W\in M_n(\mathbb{R})\), the idetity
\[
    \frac 1{|S_{n-1}|} \int_{S_{n-1}}  (y^T)^\top W y^T dS(y) = \frac{\operatorname{tr}  W}n.
\]
holds, where $S^{n-1}$ is the unit sphere in $\mathbb R^n$ and $|S^{n-1}|$ denotes its surface measure. 
Consequently, the trace of the inverted Gramian \(M^{-1}\) is proportional to the average energy required to steer the system towards terminal states on the unit sphere. 

In this section, we compare the above introduced geometric and energetic measures  for the classical and the twofold control problem. In order to make the comparison meaningful, we normalise the cost functional \(J\) by fixing \(\beta=1\). Moreover, whenever the inverse of \( G\) is considered, we assume that \( G\succ0\).

For the classical control problem, the initial state is fixed at zero, and the reachable set under the unit-energy constraint is
\begin{align}\label{eq:E_cl}
\begin{aligned}
    \mathcal E_{\mathrm{cl}} :
    & = \left\{ y^T \in \mathbb{R}^n : 
          y^T = L(0,u), \; \|(u)\|_{L^2}^2 \leq 1 \right\}\\
       & = \left\{ y^T\in\mathbb R^n \;\middle|\; (y^T)^\top  G^{-1}y^T\leq 1 \right\},
\end{aligned}
\end{align}
which represents the set of terminal states reachable from the zero initial state using distributed control inputs with energy not exceeding one and provided that $G\succ0$.

Since every admissible classical control $u$ corresponds to the admissible pair $(0,u)$ in the twofold formulation, we immediately obtain \(\mathcal E_{\mathrm{cl}} \subseteq \mathcal E_{\mathrm{ext}}\) and consequently,  \(\operatorname{vol}\!\left(\mathcal E_{\mathrm{ext}}\right) > \operatorname{vol}\!\left(\mathcal E_{\mathrm{cl}}\right)\). This shows that incorporating the optimally selected initial state always enlarges the reachable set.

This enlargement can be quantified explicitly by
the determinant ratio
\begin{align}\label{eq:ratio_vol}
    R_{\text{vol}}:= \sqrt{\frac{\det M}{\det G}}>1.
\end{align}
The ratio is well defined whenever \(G\) is positive definite. 
Since $E$ is positive definite, we have $M=\frac 1\alpha E + G\succ G \succeq 0$. Hence, by the strict monotonicity of the determinant on the cone of positive semidefinite matrices, $\sqrt{\det M}> \sqrt{\det G}$, even when G is singular.

Besides enlarging the reachable set, the twofold controllability Gramian also offers an energetic advantage. The following theorem establishes that the average minimum cost of reaching terminal states on the unit sphere is strictly lower than in the classical formulation.

\begin{teo}
    For every  operator pair \((A,B)\) such that \( G\succ0\),
    \begin{align}\label{eq:ratio_tr_inv}
        R_{\text{tr}}: =\frac{\operatorname{tr}\!\left(M^{-1}\right)}{\operatorname{tr}\!\left( G^{-1}\right)} < 1.
    \end{align}
\end{teo}

\begin{proof}
    Since \(M \succ  G \succ 0\), the order-reversing property of matrix inversion gives \(M^{-1} \prec G^{-1}\).
    
    Hence, \( G^{-1} - M^{-1} \succ 0\). Taking the trace yields
    \[
        \operatorname{tr} \! \left(M^{-1} \right) < \operatorname{tr} \! \left( G^{-1} \right).
    \]
\end{proof}

\begin{nap}
In the case of a singular Gramian $ G $, the average control cost (required to reach a target on the unit sphere) in the classical formulation  is infinite due to the presence of unreachable directions. In this case the last statement still remains valid if $\operatorname{tr}\! \left(G^{-1}\right)$ is rewritten as  $\sum_i \left(\sigma_i( G)\right)^{-1}$, where $\sigma_i( G)$ denote the eigenvalues of the Gramian, and the inverse of zero is treated as infinity.
\end{nap}

\section{Convergence Analysis of the Twofold Control as the Cost of Adjusting Initial Datum Goes to Infinity}

In this section, we analyse the convergence of the optimal control cost as the weight \(\alpha\) associated with the initial state tends to infinity, while keeping  \(\beta = 1\) fixed. The aim is to investigate whether increasing the cost of adjusting the initial datum encourages the system to start from zero and rely exclusively on the distributed control.

The associated cost functional is given by
\begin{align*}
    J_{\alpha}(y^0, u) =  \alpha \|y^0\|_{\mathbb{R}^n}^2 + \|u\|_{L^2(0,T;\mathbb{R}^m)}^2 , \qquad \alpha > 0.
\end{align*}
Since \(J_\alpha\) is obtained from \eqref{eq:cost_ab} by fixing \(\beta = 1\), the existence and uniqueness result established for \(J\) immediately implies the existence and uniqueness of a minimiser of \(J_\alpha\). We therefore seek the minimising pair
\[
    (\tilde y_{\alpha}^0, \tilde u_{\alpha}) = \operatorname*{arg\,min}_{(y^0,u) \in X} J_{\alpha}(y^0, u),
\]
subject to the final state constraint \eqref{eq:Lu=xT}.

In the same manner, we denote the corresponding twofold controllability Gramian
\begin{align}\label{eq:Malpha}
    M_\alpha := \frac{1}{\alpha} E +  G.
\end{align}

To investigate the asymptotic behaviour of the optimal controls as \(\alpha \to \infty\), we first analyse the properties of the operator \(M_\alpha\) appearing in \eqref{eq:Malpha}. We distinguish two cases: the one with a positive definite controllability Gramian \( G\), and a general one with no additional assumptions on \( G\).

\subsection{Case: The Gramian Operator is Positive Definite}

Firstly, we consider the case in which the controllability Gramian \( G\) is positive definite. Since the first term of \eqref{eq:Malpha} is bounded and converges to zero in operator norm as \(\alpha \to \infty\), we have
\[
    M_\alpha \to M_\infty :=  G
\]
in the same norm. Moreover, as \( G\) is positive definite and the inversion map is continuous on the set of invertible matrices, we have
\[
    M_\alpha^{-1} \to M_\infty^{-1} :=  G^{-1}
\]
in operator norm.

This convergence allows us to characterise the asymptotic behaviour of the optimal adjoint state, the optimal  control inputs, and the optimal cost function.

\begin{teo}\label{thm-conv}
    Let \(\tilde p_\alpha^T = M_\alpha^{-1} y^T\), where \(y^T \in \mathbb R^n\) is a fixed final state, and let $(\tilde  y_\alpha^0, \tilde  u_\alpha)$ be the corresponding optimal control pair given by \eqref{eq:(y0,u)=L*pT}, i.e. 
    \begin{align}\label{eq:(y_0,u)_alpha}
           \tilde  y_\alpha^0 := \frac{1}{\alpha} e^{T A^\star} \tilde  p_\alpha^T, \qquad \tilde u_\alpha(t) := B^\star e^{(T-t) A^*} \tilde p_\alpha^T,
    \end{align}
    for \(t \in [0, T]\). 
    
    Then, as \(\alpha \to \infty\), the following convergence results hold:
    \begin{enumerate}
        \item \(\tilde  p_\alpha^T \to \tilde  p_\infty^T: = M_\infty^{-1} y^T\quad \text{%strongly 
        in } \mathbb R^n.\)
        \item \(\tilde  y_\alpha^0 \to 0 \qquad \text{ in } \mathbb R^n.\)
        \item \(\tilde  u_\alpha \to \tilde  u_\infty:= B^* e^{(T-t) A^* } \tilde p_\infty^T \quad \text{strongly in } L^2([0, T]; \mathbb R^m),\)
        and  $\tilde  u_\infty$ is the classical optimal control steering the system from zero to $y^T$, i.e. \(L(0, \tilde  u_\infty)=y^T\). 
        \item The corresponding costs satisfy
        \[
            \lim_{\alpha \to \infty} J_\alpha(\tilde  y_\alpha^0, \tilde  u_\alpha) = J_\alpha(0, \tilde u_\infty)=   \| \tilde u_\infty \|_{L^2(0,T;\mathbb{R}^m)}^2,
        \]
        where the convergence is from below. 
    \end{enumerate}
\end{teo}

\begin{proof}
    Since \(M_\alpha^{-1} \to M_\infty^{-1} =  G^{-1}\) is in operator norm, for the fixed terminal state \(y^T \in \mathbb R^n\) we obtain
    \[
        \tilde p_\alpha^T = M_\alpha^{-1} y^T \to M_\infty^{-1} y^T = \tilde p_\infty^T \quad \text{in } \mathbb R^n.
    \]
    Since \(e^{T A^\star}\) is a bounded operator and \((\tilde p_\alpha^T)_\alpha\) is convergent, hence bounded, it follows that the optimal initial datum given by \eqref{eq:(y_0,u)_alpha} satisfies 
    \[
        |\tilde y_\alpha^0|_{\mathbb R^n} \le \frac{C}{\alpha} |\tilde p_\alpha^T|_{\mathbb R^n} \to 0,
    \]
    which provides the second announced claim.
    
    Moreover,
    \[
        \tilde u_\alpha(t) - \tilde u_\infty (t) = B^\star e^{(T-t)A^\star} (\tilde p_\alpha^T - \tilde p_\infty^T).
    \]
    Since the operator \(B^\star e^{(T-t)A^\star}\) is bounded from \(\mathbb R^n\) into \(L^2(0, T; \mathbb R^m)\), the convergence \(\tilde p_\alpha^T \to \tilde p_\infty^T\) implies
    \[
        \tilde u_\alpha \to \tilde u_\infty \qquad \text{strongly in } L^2(0, T; \mathbb R^m).
    \]
    By plugging the relation $\tilde  p_\infty^T: = M_\infty^{-1} y^T$ into the definition of $\tilde u_\infty$, we get 
    \[
    \int_0^T e^{(T - t)A} B \tilde u_\infty (t) \, dt =M_\infty \tilde p_\infty^T = y^T,
    \]
    which shows that $\tilde  u_\infty$ is the classical optimal control steering the system from zero to $y^T$.
    
    Finally, the convergence of the cost functional 4) follows directly from the previous convergences and from
    \[
        \alpha \|\tilde y_\alpha^0\|^2 = \frac{1}{\alpha} \|e^{T A^\star} \tilde p_\alpha^T\|^2 \to 0.
    \]
    As 
    \[
        J_\alpha(\tilde y_\alpha^0, \tilde u_\alpha) =\min_{(y^0, u) \in \mathcal A} J_{\alpha}( y^0, u) \leq  J_{\alpha}(0, \tilde u_\infty)=  \|\tilde u_\infty\|^2 ,
    \]
    the optimal cost values $J_\alpha(\tilde y_\alpha^0, \tilde u_\alpha)$ converge from below to its limit value. 
\end{proof}

Theorem \ref{thm-conv} reveals the twofold control formulation %can be regarded 
as a genuine extension of the classical minimum-energy control problem. As the penalty on the initial state tends to infinity, the contribution of the optimised initial condition vanishes, and the optimal solution converges to the classical minimum-energy control $\tilde u_\infty$ associated with the controllability Gramian \( G\) and the prescribed initial condition \(y^0 = 0\).

\subsection{Case: The Gramian Operator is Positive Semidefinite}

We now consider the case in which the controllability Gramian \( G\) is a positive semidefinite matrix. Unlike in the positive definite case, \( G\) may have a non-trivial kernel, and the asymptotic behaviour of the optimal solution depends on the position of the target state \(y^T\) relative to \(\operatorname{Im}( G)\) and \(\ker( G)\). To exploit this structure, we use the orthogonal decomposition
\begin{align}\label{eq:ort_decomp}
    \mathbb{R}^n  =\overline{\operatorname{Im}( G)} \oplus \ker( G).
\end{align}
Using the singular value decomposition \cite[Theorem 7.3.2]{HJ13}, one can write the Gramian as 
\[
     G = USU^\star,
\]
where $U$ is an orthonormal matrix whose first $r\leq n$ columns of \(U\) form an orthonormal basis of \(\operatorname{Im}( G)\), while the remaining ones form an orthonormal basis of \(\ker( G)\). Consequently,
\begin{align}\label{eq:decomp_G}
    U^\star  G U = S:=
    \begin{bmatrix}
        S_r & 0\\ 0 & 0
    \end{bmatrix},
\end{align}
where \(S_r = \operatorname{diag} (\sigma_1, \ldots, \sigma_r)\) is a positive definite diagonal matrix containing the non-zero singular values of \(G\). In addition, we shall decompose $U$ as 
    \begin{align}
        \label{eq:U_decomp}
            U=[U_1\ U_2],
   \end{align}
    where the columns of \(U_1\) and \(U_2\) form orthonormal bases of \(\operatorname{Im}( G)\) and \(\ker( G)\), respectively.

The positive definite operator \(E = e^{TA}e^{TA^\star}\) is represented in the same basis as
\begin{align}\label{eq:decomp_E}
    U^\star E U =
    \begin{bmatrix}
        E_{11} & E_{12}\\ E_{21} & E_{22}
    \end{bmatrix}.
\end{align}
The positivity of \(E\) implies that \(E_{22}\) is positive definite, and the Schur complement \(E_{11} - E_{12}E_{22}^{-1}E_{21}\) is also positive definite \cite[Ch. 6]{ZF05}. 

The representation of the target state in the orthogonal basis induced by \(U\) is
\begin{align}\label{eq:yhat}
    \hat y^T:=U^\star y^T=
    \begin{bmatrix}
        \hat y_1^T\\ \hat y_2^T
    \end{bmatrix},
\end{align}
where \(\hat y_1^T\in\mathbb{R}^r\) and \(\hat y_2^T\in\mathbb{R}^{n-r}\) denote coefficients corresponding to basis vectors  from \(\operatorname{Im}( G)\) and \(\ker( G)\), respectively. If \(y^T\in\operatorname{Im}( G)\), then \(\hat y_2^T=0\).

Similarly, we represent the adjoint state in the same orthogonal basis,
\[
    \hat p_\alpha^T:=U^\star \tilde p_\alpha^T= \begin{bmatrix} \hat p_{\alpha,1}^T\\ \hat p_{\alpha,2}^T \end{bmatrix},
\]
where the first block corresponds to \(\operatorname{Im}( G)\) and the second to \(\ker( G)\).

We can now state the next theorem that provides the desired  convergence result.

\begin{teo}\label{thm-conv_singG}
    Let \(M_{\alpha}\) be defined by \eqref{eq:Malpha} and  \(\tilde p_\alpha^T = M_\alpha^{-1} y^T\). Then, as \(\alpha\to\infty\), the following statements hold.
    \begin{enumerate}
        \item If \(y^T\in\operatorname{Im}( G)\), then %there exists \(p_\infty^T\in\mathbb R^n\) such that
        \[
            \tilde p_\alpha^T\longrightarrow \tilde p_\infty^T: =U \hat p_\infty^T.
        \]
       where
        \[
            \hat p_\infty^T=
            \begin{bmatrix}
                S_r^{-1}\hat y_1^T\\[1mm] -E_{22}^{-1}E_{21}S_r^{-1} \hat y_1^T
            \end{bmatrix},
        \]
        with \(S_r\), \(E_{22}\), and \(E_{21}\) the blocks given by the decompositions \eqref{eq:decomp_G} - \eqref{eq:decomp_E} and $y_1^T$ given by \eqref{eq:yhat}.
        
        Furthermore,
        \[
            \tilde y_\alpha^0= \frac1\alpha e^{TA^\star}\tilde p_\alpha^T \longrightarrow0,
        \]
        and
        \[
            \tilde u_\alpha(t)\longrightarrow \tilde u_\infty(t):= B^\star e^{(T-t)A^\star}\tilde p_\infty^T \quad\text{in }L^2(0,T;\mathbb{R}^m),
        \]
        and  $\tilde  u_\infty$ is the classical optimal control steering the system from zero to $y^T$, i.e. \(L(0, \tilde  u_\infty)=y^T\). 
        
        Finally,
        \[
            J_\alpha(\tilde y_\alpha^0, \tilde u_\alpha) \longrightarrow J(0,\tilde u_\infty).
        \]
        \item If \(y^T\notin\operatorname{Im}(G)\), let
        \[
            y^T=y_1^T+y_2^T, \qquad y_1^T\in\operatorname{Im}( G), \quad 0\neq y_2^T\in\ker( G).
        \]
        Then, as \(\alpha\to\infty\),
        \begin{equation}
            \label{conv:p1}
                     \hat p_{\alpha,1}^T \longrightarrow \hat p_{\infty,1}^T := S_r^{-1} \left( \hat y_1^T-E_{12}E_{22}^{-1}\hat y_2^T \right),  
      \end{equation}
        whereas
       \begin{equation}
            \label{conv:p2}
        \|\hat p_{\alpha,2}^T\|\longrightarrow\infty, \qquad \|\widetilde p_\alpha^T\|\longrightarrow\infty
       \end{equation}
        Furthermore, the optimal initial states satisfy
        \begin{align}
        \label{eq:conv_y0}
                    \tilde y_\alpha^0\to \tilde y_\infty:= e^{-TA}(y_2^T+h),
               \end{align}
                      with
        \begin{align}\label{eq:corection-h}
            h:=U_1E_{12}E_{22}^{-1}\hat y_2^T \in\operatorname{Im}( G),
        \end{align}
        while for the optimal controls it holds
        \begin{equation}
            \label{conv:u}
            \tilde u_\alpha\to \tilde u_\infty := B^\star e^{(T-t)A^\star} U_1 \hat p_{\infty,1}^T\quad\text{in }L^2(0,T;\mathbb R^m).     \end{equation}      
        Moreover,  \(\tilde u_\infty\) is the classical minimum-energy control steering the system from zero to \(y_1^T-h\), i.e.,
        \(
            L(0,\tilde u_\infty)=y_1^T-h.
        \)
        Consequently, the limit control pair $(\tilde y_\infty^0,\tilde u_\infty)$ also steers the system to the prescribed target $y^T$, while 
        \[
            J_\alpha(\tilde y_\alpha^0,\tilde u_\alpha)\to\infty.
        \]
    \end{enumerate}
\end{teo}

\begin{proof}
    Let \(y^T\in \mathbb{R}^n\) be a given target state, and let \(\hat y^T\) be its representation in the orthogonal basis defined in \eqref{eq:yhat}. The corresponding block representations of \( G\) and \(E\) are given by \eqref{eq:decomp_G} and \eqref{eq:decomp_E}, respectively,  while the equation \eqref{eq:xT=MpT} for the optimal adjoint state can be rewritten in block form as
    \begin{align}\label{eq:blockform}
        \begin{cases}
            \left(S_r + \frac{1}{\alpha} E_{11}\right)\hat p_{\alpha,1}^T + \frac{1}{\alpha} E_{12}\hat p_{\alpha,2}^T = \hat y^T_1,\\[1mm]
            \frac{1}{\alpha} E_{21}\hat p_{\alpha,1}^T + \frac{1}{\alpha} E_{22}\hat p_{\alpha,2}^T = \hat y^T_2.
        \end{cases}
    \end{align}
    
    1) If \(y^T \in \operatorname{Im}( G)\), then the target state has no component in \(\ker( G)\), i.e. $\hat y^T_2=0$.         
    Substituting \(\hat p_{\alpha,2}^T\) from the second equation of the  block system \eqref{eq:blockform} into the first yields
    \[
        \left(S_r + \frac{1}{\alpha}(E_{11} - E_{12}E_{22}^{-1}E_{21})\right)\hat p_{\alpha,1}^T = \hat y_1^T.
    \]
    Since \(S_r\) and the Schur complement \(E_{11} - E_{12}E_{22}^{-1}E_{21}\) are positive definite, the above matrix is invertible for every \({\alpha}  > 0\). Therefore,
    \[
        \hat p_{\alpha,1}^T = \left(S_r + \frac{1}{\alpha}(E_{11} - E_{12}E_{22}^{-1}E_{21})\right)^{-1} \hat y_1^T.
    \]
    Passing to the limit as \({\alpha} \to \infty\) and using the continuity of matrix inversion, we obtain
    \[
        \hat p_{\alpha,1}^T \to S_r^{-1}\hat y_1^T.
    \]
    Consequently,
    \[
        \hat p_{\alpha,2}^T = -E_{22}^{-1}E_{21}\hat p_{\alpha,1}^T \to - E_{22}^{-1}E_{21}S_r^{-1}\hat y_1^T.
    \]
    Hence,
    \[
        \widetilde p_\alpha^T \to \widetilde p_\infty^T = U
        \begin{bmatrix}
            S_r^{-1}\hat y_1^T\\ -E_{22}^{-1}E_{21}S_r^{-1}\hat y_1^T
        \end{bmatrix}.
    \]
    The convergence of \(y_\alpha^0\), \(u_\alpha\), and \(J_\alpha\) then follows exactly as in the positive definite case (Theorem \ref{thm-conv}).

    2) Suppose now that \(y^T \notin \operatorname{Im}( G)\). Then the target state has a nonzero component in \(\ker( G)\), that is, \(\hat y_2^T\neq 0\).
   
    Since \(E\succ0\), its principal submatrix \(E_{22}\) is positive definite and hence invertible. From the second equation of the  block system \eqref{eq:blockform}, we obtain
    \begin{align}\label{eq:palpha2}
        \hat p_{\alpha,2}^T = E_{22}^{-1} \left( \alpha \hat y_2^T -E_{21}\hat p_{\alpha,1}^T\right).
    \end{align}
    Substituting this expression into the first equation gives
    \begin{align*}
        \hat p_{\alpha,1}^T = \left( S_r+\frac{1}{\alpha} \left(E_{11}-E_{12}E_{22}^{-1}E_{21}\right) \right)^{-1} \left( \hat y_1^T-E_{12}E_{22}^{-1}\hat y_2^T \right).
    \end{align*}
    Passing to the limit as \(\alpha\to\infty\), we obtain
   \eqref{conv:p1}.
   In particular,  the component of the optimal adjoint datum belonging to \(\operatorname{Im}( G)\) remains bounded.
    
    On the other hand, dividing the expression \eqref{eq:palpha2} by \(\alpha\) and using the boundedness of $(\hat p_{\alpha,1}^T)$ yields
    \begin{align}\label{eq:palpha2limit}
        \frac{1}{\alpha}\hat p_{\alpha,2}^T = E_{22}^{-1}\hat y_2^T - \frac{1}{\alpha}E_{22}^{-1}E_{21} \hat p_{\alpha,1}^T \longrightarrow E_{22}^{-1}\hat y_2^T.
    \end{align}
    Since \(\hat y_2^T\neq0\) and \(E_{22}\) is invertible, \(E_{22}^{-1}\hat y_2^T\neq0\). Consequently, the divergence of the optimal adjoint data \eqref{conv:p2} follows.

   In order to establish the convergence of the optimal initial states, we combine \eqref{conv:p1} and \eqref{eq:palpha2limit}, thus obtaining
    \[
        \frac{1}{\alpha}\tilde p_\alpha^T = U
        \begin{bmatrix}
            \frac{1}{\alpha}\hat p_{\alpha,1}^T\\[1mm]
            \frac{1}{\alpha}\hat p_{\alpha,2}^T
        \end{bmatrix}
        \longrightarrow U
        \begin{bmatrix}
            0\\
            E_{22}^{-1}\hat y_2^T
        \end{bmatrix}.
    \]
    Consequently, the optimal initial state satisfies
    \begin{align}
        \label{eq:y_inf}
        \tilde y_\alpha^0 = \frac{1}{\alpha} e^{TA^\star}\tilde p_\alpha^T \longrightarrow \tilde y_\infty^0 := e^{TA^\star}U
        \begin{bmatrix}
            0\\
            E_{22}^{-1}\hat y_2^T
        \end{bmatrix}
        \not= 0.
    \end{align}
   Using the decomposition \eqref{eq:U_decomp} it follows 
    \begin{align}\label{eq:eTAy0inf}
        \begin{aligned}
            \tilde y_\infty^0 &= e^{TA^\star}U_2E_{22}^{-1}\hat y_2^T\\
            &= e^{-TA}(U_1 U_1^\star + U_2 U_2^\star )E U_2E_{22}^{-1}\hat y_2^T  \\
            &= e^{-TA}(U_1E_{12}E_{22}^{-1}\hat y_2^T + U_2\hat y_2^T) \\
            & = e^{-TA}(h + y_2^T),
        \end{aligned}
    \end{align}
     with $h$ given by \eqref{eq:corection-h}.
  
   Finally, since \(U_2\hat p_{\alpha,2}^T \in\ker( G)\), for the optimal distributed controls we have
    \[
    \begin{aligned}
        \tilde u_\alpha(t)
        &= B^\star e^{(T-t)A^\star} \tilde p_\alpha^T\\
        &= B^\star e^{(T-t)A^\star} \left( U_1\hat p_{\alpha,1}^T + U_2\hat p_{\alpha,2}^T \right)\\
        &= B^\star e^{(T-t)A^\star} U_1\hat p_{\alpha,1}^T.
    \end{aligned}
    \]
    The asserted convergence \eqref{conv:u} now follows directly from \eqref{conv:p1}. 
   
    Using the derived expression for $\tilde u_\infty$ we get
       \begin{align}\label{eq:L(0,uinf)}
        \begin{aligned}
            L(0,\tilde u_\infty) &=  GU_1S_r^{-1}(\hat y_1^T-E_{12}E_{22}^{-1}\hat y_2^T) \\
            &= U_1 \hat y_1^T- U_1E_{12}E_{22}^{-1}\hat y_2^T \\
            &:=  y_1^T - h,
        \end{aligned}
    \end{align}
    showing that  \(\tilde u_\infty\) is the unique minimum-energy control steering the system from zero to \(y_1^T-h\).
    
    Combining \eqref{eq:eTAy0inf} and \eqref{eq:L(0,uinf)}, we obtain
    \[
        \begin{aligned}
            L(\tilde y_\infty^0,\tilde u_\infty)
            &=e^{TA}\tilde y_\infty^0+L(0,\tilde u_\infty)\\
            &=(h+y_2^T)+(y_1^T-h)\\
            &=y_1^T+y_2^T=y^T.
        \end{aligned}
    \]
    Thus, the limiting pair satisfies the terminal-state constraint \eqref{eq:Lu=xT}.

   Consequently, convergences \eqref{conv:u} and \eqref{eq:y_inf} imply
    \[
        J_\alpha(\tilde y_\alpha^0,\tilde u_\alpha)= \alpha\|\tilde y_\alpha^0\|^2+\|\tilde u_\alpha\|_{L^2}^2 \longrightarrow+\infty.
    \]

    Therefore, although both $\|\tilde y_\alpha^0\|^2$ and $\|\tilde u_\alpha\|^2$ converge to finite values, the corresponding control costs diverge as $\alpha\to\infty$.
\end{proof}

The last theorem shows that if the target \(y^T\) is reachable by a distributed control only, i.e. if it belongs to \(\operatorname{Im}( G)\), then the optimal control pair satisfies \((\tilde y_\alpha^0, \tilde u_\alpha)  \mathop{ \longrightarrow }_{\alpha} (0, \tilde u_\infty)\), where \(\tilde u_\infty\) is the classical optimal control steering the system from zero to  \(y^T\). In other words, the increased penalisation of adjusting the initial datum stimulates the system to choose the trivial initial datum and rely exclusively on the distributed control. 

If this is not the case, i.e. if \(y^T \not\in\operatorname{Im}( G)\), then the initial datum differs from zero and the control cost explodes as the penalisation parameter \(\alpha\) goes to infinity.

The theorem also reveals interesting and, at first glance, surprising features in the case $y^T\in\operatorname{ker}( G)$. 
Although such a target is completely unreachable in the classical formulation, the distributed control input in general remains active. It shifts the target in a way that minimizes the cost of reaching the shifted target from the initial datum alone. 

The exception is the case  \(E_{12}=0\), i.e. when the dynamics of $A$ does not mix the image and kernel of the gramian $G$.  Then the shift \(h\), given by \eqref{eq:corection-h} equals  zero, so that \(\tilde y_\infty^0=e^{-TA}y_2^T\) and \(\tilde u_\infty=0\). This behaviour will also be illustrated by the numerical example in the next section.

\section{Numerical Experiments}

In this Section we verify the above obtained theoretical results on several particular systems. In particular, we construct the optimal twofold controls, compare their cost with the classical control one, investigate the geometric and energetic controllability properties associated with the twofold controllability Gramian, and examine the asymptotic properties as the weight $\alpha$ converges to $\infty$.  

\subsection{Case: The Gramian Operator is Positive Definite}

\subsubsection{1D Heat Equation With a Boundary Control}

We first apply the theory to a semi-discretised one-dimensional heat equation with boundary control on the right-hand side. To this aim, we consider the system \eqref{eq:system}, where the operator \(A\) and control matrix \(B\) are given by
\[
    A = (n+1)^2\,\mathrm{tridiag}(1,-2,1), \quad  B = (n+1)^2 e_n,
\]
where \(e_n\) is the \(n\)-th canonical basis vector of \(\mathbb{R}^n\).

This corresponds to a finite-difference spatial discretisation of the one-dimensional heat equation with \(n\) nodes, where control is applied at the right boundary \((x = L)\)
\[
    \left\{
    \begin{aligned}
        &\frac{\partial y}{\partial t}(x,t) = \frac{\partial^2 y}{\partial x^2}(x,t), \quad & x \in (0,L),\, t > 0, \\
        &y(0,t) = 0,\quad y(L,t) = u(t), & t \geq 0, \\
        &y(x,0) = y^0(x),\quad & x \in [0,L].
    \end{aligned}
    \right.
\]

We set the time horizon \(T = 0.1\), the interval length \(L = 1\), and we use a uniform grid of \(n = 8\) nodes. The both control cost parameters \(\alpha\) and \(\beta\) are set to \(1\). After discretising the heat equation,  the obtained system is controllable in the classical setting, and we can steer it to the desired final state, which in this example is chosen as the (discretised) sine function, 
\[
    y_i^T=\sin\!\left(\dfrac{\pi i}{n+1}\right), \qquad  i=1,\ldots,n.
\]

This enables us to compare two approaches: the standard distributed control one and the twofold control. In particular, we aim to compare the corresponding optimal control costs and check to which level does the introduction of the latter notion reduces the effort required to reach the target. 

\begin{figure}[htbp]
    \centerline{\includegraphics[width=\columnwidth]{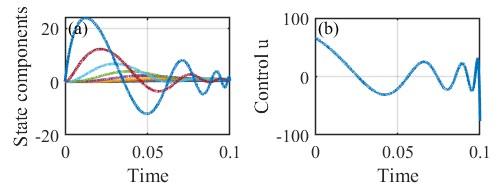}}
    \caption{Evolution of the optimal state and control function for the heat equation with the distributed control input only.}
    \label{fig1}
\end{figure}

Fig. \ref{fig1} shows the evolution of the system states and the control function over time in the classical control setting, where the system is influenced through the distributed control only, while the initial datum is fixed to zero.  For this approach, the optimal value of the cost functional is \(72\).

\begin{figure}[htbp]
    \centerline{\includegraphics[width=\columnwidth]{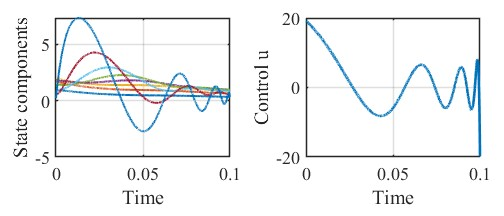}}
    \caption{Optimal solution for the heat equation obtained by the proposed twofold control strategy: (a) evolution of the state components and (b) corresponding distributed control input.}
    \label{fig2}
\end{figure}

Fig. \ref{fig2} shows the same scenario, but in this case, the system is controlled by both the initial state and the distributed control. The corresponding optimal value of the cost functional is \(20.282\). When we compare the optimal costs  for distributed control with those for twofold control, we observe a substantial reduction – it is 3.5
times lower.

To further assess the benefits of the twofold control approach, we compare the corresponding controllability Gramians using the volumetric and average energy measures. For the heat equation, the numerical results yield \(R_{\text{vol}}\approx 1.904 > 1\) and \(R_{\text{tr}}\approx 0.685 < 1\), in agreement with the theoretical inequalities \eqref{eq:ratio_vol} and \eqref{eq:ratio_tr_inv}, respectively. These results confirm that the twofold control framework enlarges the reachability ellipsoid while reducing the average energy required to reach a target state, providing further evidence of its advantages over the classical control setting.

In the next step we analyse the dependence of the optimal cost on the weighting  parameter \(\alpha\).

To this aim, the values of the parameter were selected from a logarithmically equally spaced sequence of \(20\) values between \(10^0\) and \(10^8\).

\begin{table}[H]
    \caption{Optimal initial-state and control energies and the corresponding cost functional values for the heat equation.}
    \label{table1}
    \centering
        \begin{tabular}{|c|c|c|c|c|} \hline
            $i$ & $|\tilde y^0|^2_{\mathbb{R}^n}$ & $\|\tilde u\|^2_{L^2}$ & $J(\tilde y^0,\tilde u)$ \\ \hline
            $1$ & $14.5$ & $5.7$ & $20.3$ \\ \hline
            $3$ & $2.0$ & $38.5$ & $52.6$ \\ \hline
            $5$ & $0.1$ & $64.9$ & $68.4$ \\ \hline
            $10$ & $5\cdot 10^{-6}$ & $71.9$ & $72.0$ \\ \hline
            $20$ & $2\cdot 10^{-14}$ & $72.0$ & $72.0$ \\ \hline
        \end{tabular}
\end{table}

\begin{figure}[htbp]
    \centerline{\includegraphics[width=\columnwidth]{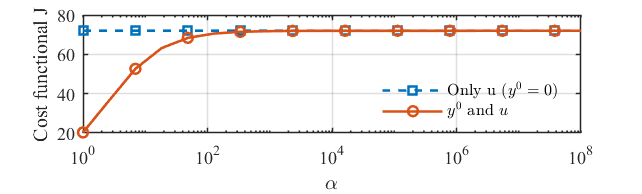}}
    \caption{Cost functional for the classical and twofold control problems as a function of the penalty parameter \(\alpha\).}
    \label{fig3}
\end{figure}

Table \ref{table1} and Fig. \ref{fig3} show convergence results with respect to  increasing \(\alpha\).
%for both the heat and wave equations. 
Table \ref{table1} summarises the norms of the optimal initial state, the control input, and the total cost, where \(i\) denotes the index of the corresponding value of \(\alpha\) in the considered sequence. In Fig. \ref{fig3}, markers indicate the values corresponding to every second value of \(\alpha\) in this sequence.

Note that by taking large values of \(\alpha\), the initial state tends to zero, and the distributed control input dominates. Although theoretically, the optimal initial value for any $\alpha\in \mathbb{R}^+$ is non-zero, we note that already for \(\alpha=10^8\), its norm is at the level of machine precision. Consequently, the numerical solution in that case is practically indistinguishable from the classical control solution with \(y^0=0\). Such behaviour is expected, as the high cost of adjusting the initial state fosters the system to start from zero. 

\subsubsection{1D Wave Equation With a Boundary Control}

In this example we consider the system arising from a finite-difference semi-discretisation of the one-dimensional wave equation. Let \(n\) denote the number of interior spatial grid points. Rewriting the second-order wave equation as a first-order system yields a state vector in \(\mathbb{R}^{2n}\), whose first \(n\) components represent the displacement and the remaining \(n\) components represent the velocity. The corresponding matrices entering the system \eqref{eq:system} are
\[
    A = 
    \begin{bmatrix}
         0      & I \\
          \, \Delta & 0
    \end{bmatrix}, \quad B = (n+1)^2 e_{2n},
\]
where \(e_{2n}\) is the \(2n\)-th canonical basis vector in \(\mathbb{R}^{2n}\), corresponding to the displacement at the last interior node, and \(\Delta\) is the standard finite-difference Laplacian, \(\Delta = (n+1)^2 \, \mathrm{tridiag}(1,-2,1)\).

This setting corresponds to the one-dimensional wave equation with boundary control applied at the right endpoint \((x = L)\)
\[
    \left\{
    \begin{aligned}
        &\frac{\partial^2 y}{\partial t^2}(x,t) = \frac{\partial^2 y}{\partial x^2}(x,t), \quad && x \in (0,L),\; t > 0, \\
        &y(0,t) = 0,\quad y(L,t) = u(t), && t \geq 0, \\
        &y(x,0) = y^0(x),\quad \frac{\partial y}{\partial t}(x,0) = v^0(x), && x \in [0,L].
    \end{aligned}
    \right.
\]

The parameter setting for this experiment is chosen as \(T = 1\), \(L = 1\), \(n = 4\), and \(\alpha = 1\). The final state is prescribed as 
\[
    y_i^T=
    \begin{cases}
        \sin\!\left(\dfrac{\pi i}{n+1}\right), & i=1,\ldots,n,\\[2mm]
        0, & i=n+1,\ldots,2n.
    \end{cases}
\]
    
\begin{figure}[htbp]
    \centerline{\includegraphics[width=\columnwidth]{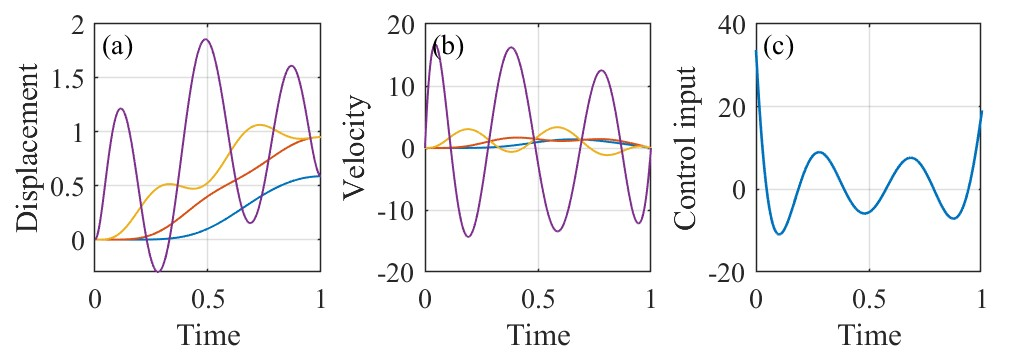}}
    \caption{Distributed control of the one-dimensional wave equation: (a) displacement components, (b) velocity components, and (c) distributed control input.}
    \label{fig4}
\end{figure}
We first fix the initial data to zero and construct the corresponding optimal distributed control (Fig. \ref{fig4}). The corresponding cost functional value equals \(52.7\).

\begin{figure}[htbp]
    \centerline{\includegraphics[width=\columnwidth]{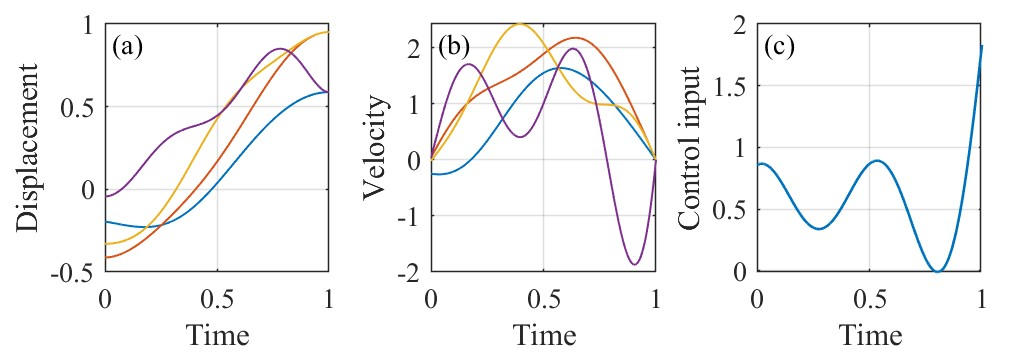}}
    \caption{Optimal twofold control of the one-dimensional wave equation for $\alpha=1$: (a) displacement components, (b) velocity components, and (c) distributed control input.}
    \label{fig5}
\end{figure}

In the next step we apply the twofold control to our system. The evolution of the optimal state and the control function is depicted in Fig. \ref{fig5}. Although they follow a similar oscillatory pattern characteristic for the wave equation as on Fig. \ref{fig4}, the amplitudes are much smaller, suggesting the significant decrease of the control cost. Indeed, the corresponding cost functional  value equals \(0.8\), which is \(66\) times smaller than in the classical control setting. 

To complement the comparison of optimal control costs, we examine the volumetric and average energy measures for the wave equation. The computed ratios are \(R_{\text{vol}}\approx 22902 > 1\) and \(R_{\text{tr}}\approx 0.0006 < 1\), respectively, which are consistent with the theoretical inequalities \eqref{eq:ratio_vol} and \eqref{eq:ratio_tr_inv}. These results further illustrate the advantages of the twofold control approach, showing an enlargement of the reachability ellipsoid together with a reduction in the average energy required to reach a target state.

\begin{table}[H]
    \caption{Optimal initial-state and control energies and the corresponding cost functional values for the wave equation.}
    \label{table2}
    \centering
        \begin{tabular}{|c|c|c|c|c|} \hline
            $i$ & $|\tilde y^0|^2_{\mathbb{R}^n}$ & $\|\tilde u\|^2_{L^2}$ & $J(\tilde y^0,\tilde u)$ \\ \hline
            $1$ & $0.38$ & $0.4$ & $0.8$ \\ \hline
            $3$ & $0.06$ & $1.2$ & $1.6$ \\ \hline
            $5$ & $0.02$ & $1.6$ & $3.0$ \\ \hline
            $10$ & $1 \cdot 10^{-3}$ & $33.0$ & $41.2$ \\ \hline
            $20$ & $9 \cdot 10^{-12}$ & $52.7$ & $52.7$ \\ \hline
        \end{tabular}
\end{table}

\begin{figure}[htbp]
    \centerline{\includegraphics[width=\columnwidth]{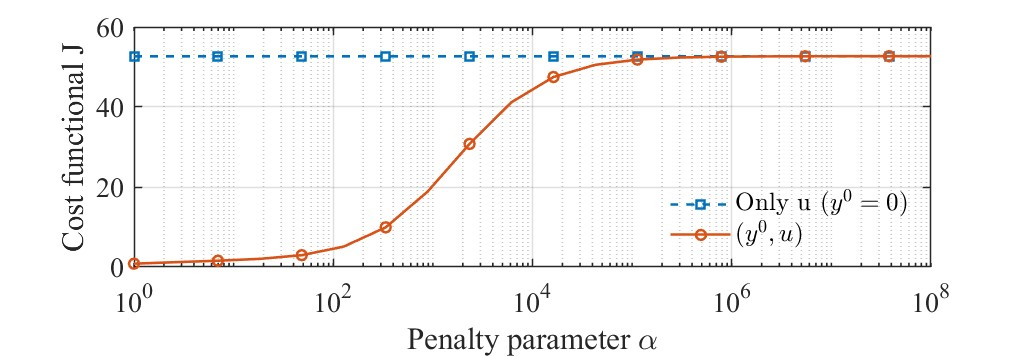}}
    \caption{Convergence of cost function for the classical and twofold control problems with respect to the parameter $\alpha$ for the wave equation.}
    \label{fig6}
\end{figure}

Table \ref{table2} and Fig. \ref{fig6} demonstrate the transition from the twofold to the classical control strategy as the penalty parameter \(\alpha\) increases. For small values of \(\alpha\), the optimisation exploits the initial state, resulting in a substantial reduction of the control effort and the total cost. In particular, for \(\alpha=1\), the optimal cost is two orders of magnitude smaller than that of the corresponding classical control problem.

As \(\alpha\) increases, the contribution of the initial state gradually diminishes, while the control energy approaches the value obtained for the classical formulation. This behaviour is clearly illustrated in Fig. \ref{fig6}, where the cost functional converges to the constant cost of the distributed control problem.

For \(\alpha=10^8\), the norm of the optimal initial state is of the order of machine precision, so the computed solution is numerically indistinguishable from the classical control solution.

\subsubsection{Comparison Between the Heat and Wave Equations}

The significant difference between the improvements obtained for the heat and wave equations can be attributed to the fundamentally different nature of their dynamics. The heat equation is strongly dissipative: due to diffusion, the influence of the initial condition decays exponentially fast, and the system gradually loses memory of its initial state. Consequently, even when the initial condition is optimised through the twofold control framework, the resulting reduction in the cost functional remains relatively moderate, as most of the effort required to reach the target state must still be provided by the distributed control input. In contrast, the wave equation is conservative and oscillatory. Its energy is preserved and continuously exchanged between kinetic and potential forms, so the initial state affects the system dynamics throughout the entire control horizon. By selecting an optimal initial condition, it is possible to exploit the natural oscillatory behaviour of the system and obtain an optimal balance between two equally efficient control inputs. 

As a result, the distributed control has substantially less work to perform, leading to a dramatic reduction in the control cost. In our numerical experiments, this effect yielded an improvement of approximately the factor of \(66\) for the wave equation, compared with only the factor of \(3\)–\(4\) for the heat equation.

A similar contrast is observed in the controllability measures: while both systems satisfy the theoretical inequalities \(R_{\mathrm{vol}}>1\) and \(R_{\mathrm{tr}}<1\), the improvement is substantially more pronounced for the wave equation compared to the heat equation. These results further illustrate how the influence of the optimised initial state on the reachability ellipsoid and the average control energy depends on the underlying system dynamics.

\subsection{Case: The Gramian Operator is Singular}

To illustrate the behaviour of the proposed optimisation problem in the case of a singular controllability Gramian, we consider the  system determined by
\[
    A= \begin{bmatrix} 1&1&1\\ 0&2&0\\ 0&0&3 \end{bmatrix},
    \qquad
    B=e_1= \begin{bmatrix} 1\\0\\0 \end{bmatrix}.
\]
The final time is set to \(T=0.1\) and the prescribed target state is
\(
    y^T=e_3
\). The weighting parameter associated with the distributed control is fixed at \(\beta=1\), whereas \(\alpha\) is varied in order to investigate the asymptotic behaviour of the optimal solution.

For this choice of \(A\) and \(B\), the control acts only on the first component of the system, and the classical controllability Gramian is given by
\[
     G = \int_0^T e^{(T-t)A}BB^\top e^{(T-t)A^\top}\,dt = \frac{e^{2T}-1}{2}\,e_1e_1^\top.
\]
Consequently, the Gramian $G$ is singular and
\[
    \operatorname{Im}( G)=\operatorname{span}\{e_1\}, \qquad \ker( G)=\operatorname{span}\{e_2,e_3\}.
\]
In particular, \(y^T\in\operatorname{ker}( G)\). One might therefore expect the optimal distributed control to vanish and the optimal control pair to be given by $(e^{-TA} y^T, 0)$. 
However, direct calculation yields that the optimal initial state has the form
\begin{align}
    \label{opt_init_alpha}
    \widetilde y_\alpha^0 = \begin{pmatrix} f(\alpha)\\ 0\\ e^{-3T} \end{pmatrix}, \quad f(\alpha):= \frac{e^{-3T}-e^{-T}}{2+\alpha(1-e^{-2T})} ,
\end{align}
while the corresponding optimal control is given by
\begin{align}
    \label{opt_u_alpha}
    \widetilde u_\alpha(t) = -\frac{\alpha(e^{2T}-1)}{\alpha(e^{2T}-1)+2e^{2T}}e^{-T-t},
\end{align}
which differs from zero for any $\alpha>0$. It shifts the terminal state away from the prescribed target in the direction of the eigenvector associated with the largest eigenvalue of  $A$, which in this example is \(
    v_{\max}=0.5e_1+e_3
\) (corresponding to $\lambda_{\max}=3$). 
This shift reduces the norm of the optimal initial state  and results in the minimal value of the cost functional. 

This shift is the most announced in the limit case when $\alpha$ goes to infinity, and it is equal to the correction term defined by \eqref{eq:corection-h}, which for this example reads 
\[
    h = \frac{1-e^{-2T}}{2}e_1.
\]
As for the optimal twofold controls, the following convergences follow directly from \eqref{opt_init_alpha} and \eqref{opt_u_alpha} 
\[
    \widetilde y_\alpha^0 \longrightarrow \tilde y_\infty^0:=e^{-3T}e_3.
\]
\[
    \widetilde u_\alpha(t) \longrightarrow \tilde u_\infty(t):=-e^{-T-t}, \qquad \alpha\to\infty.
\]
It is straightforward to check  that the obtained limit $\tilde y_\infty^0$ coincides with $e^{-TA} (y^T + h)$, as expected from relation \eqref{eq:conv_y0}. Moreover, in accordance with Theorem \ref{thm-conv_singG},  $\tilde u_\infty$ is the classical minimum-energy control steering the system from zero to \(-h\).

Thus, in the limiting regime, the initial state generates the target
together with the correction \(h\), while the control $\tilde u_\infty$ compensates
precisely for this additional target component.

\begin{figure}[htbp]
    \centering
    \includegraphics[width=\linewidth]{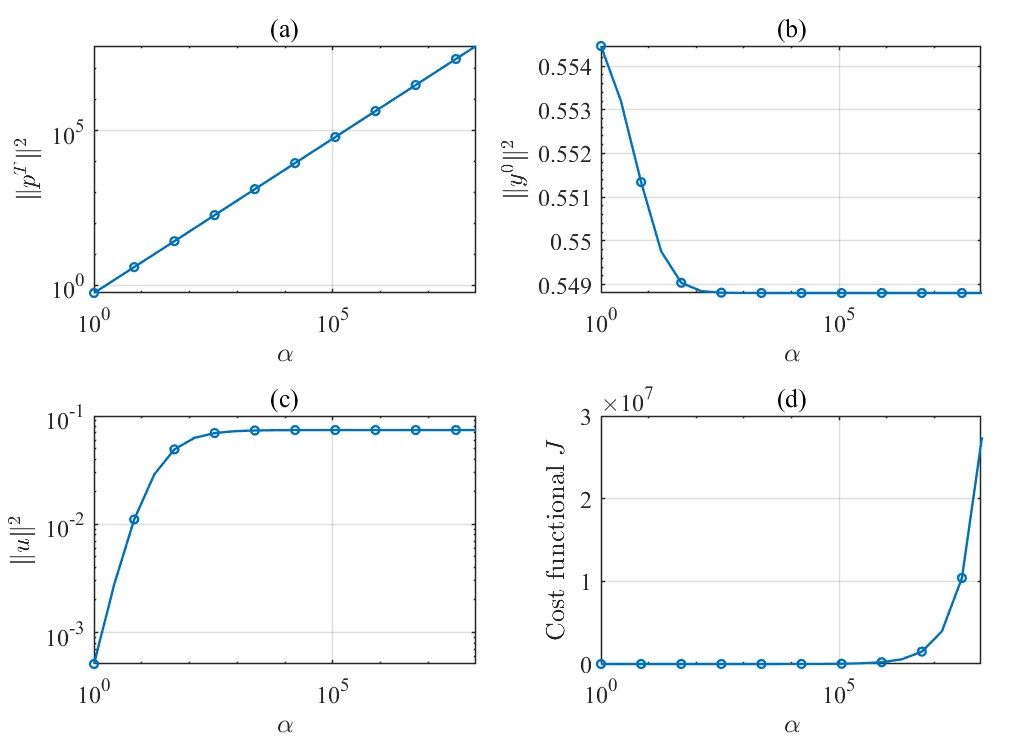}
    \caption{Dependence of the optimal quantities on the penalty parameter $\alpha$: (a) $\|p^T\|^2$, (b) $\|y^0\|^2$, (c) $\|u\|^2$, and (d) the cost functional $J$.}
    \label{fig:kvgyTnotinImG}
\end{figure}

Fig.~\ref{fig:kvgyTnotinImG} depicts the dependence of the optimal quantities on the penalty parameter $\alpha$. It confirms the asymptotic behaviour predicted by Theorem~\ref{thm-conv_singG}. As \(\alpha\to\infty\), the norm of the optimal adjoint datum \(\widetilde p_\alpha^T\) and the corresponding optimal cost diverge to infinity, whereas the optimal initial state and control remain bounded.

This example therefore clearly illustrates the contrasting asymptotic behaviours of the adjoint variable and the optimal cost, on the one hand, and the optimal initial state and control, on the other. This phenomenon occurs whenever the target state has a non-trivial component in $\ker G$, as is the case  in this example. 

\section{Conclusion}

In this paper, we introduced a twofold control framework in which the initial state and the time-dependent control input are optimised simultaneously. By allowing both mechanisms to contribute to the control process under a common cost functional, the classical optimal control problem is extended to a setting with an additional degree of freedom. Building on the HUM, we derived the corresponding optimality conditions for finite-dimensional linear systems and analysed the behaviour of the optimal pair as the penalty on the initial state tends to infinity.

A key ingredient of our analysis is the introduced twofold controllability Gramian that incorporates the contributions of both the initial state and the control input, thereby generalising the classical controllability Gramian. Unlike the latter one, the extended Gramian is always positive definite, ensuring the well-posedness of the twofold control problem and the reachability of an arbitrary target state, even when the system fails to satisfy the Kalman rank condition. In addition, comparison with the classical controllability Gramian revealed both geometric and energetic advantages of the extended formulation, including an enlargement of the reachable set under a common energy constraint and a reduction in the minimum energy required to reach a prescribed target state. Moreover, preparing a suitable initial state remains beneficial even when its cost is high: for every finite penalty, it strictly reduces the total optimal cost  compared with the classical control starting from zero.

The theoretical findings were supported by numerical experiments involving semi-discrete models of the one-dimensional heat and wave equations with boundary control. These experiments illustrated the reduction in control cost achieved by the twofold approach, demonstrated the enlargement of the reachable set, and confirmed the theoretical properties of the extended controllability framework.

Although the present analysis is restricted to finite-dimensional linear systems, the concept of twofold controllability admits natural extensions to more general dynamical settings. One important direction for future research is the investigation of infinite-dimensional systems, particularly evolution equations governed by partial differential equations. In this context, it would be of interest to establish appropriate well-posedness and reachability results and to examine whether the geometric and energetic advantages identified here persist at the continuous level. Another promising direction concerns nonlinear control systems, for which the interaction between the optimised initial state and the time-dependent control input may lead to new controllability properties and optimisation challenges. Finally, the development of efficient numerical methods for the resulting twofold optimisation problems, especially for large-scale and spatially distributed systems, constitutes a further direction for future investigation.

\end{document}